\documentclass[11pt] {article} 
\usepackage[backref]{hyperref}
  \usepackage{amsmath} 
    \usepackage{amssymb}
    
    \usepackage{pdfsync}

\newtheorem{theorem}{Theorem}[section]
\newtheorem{lemma}{Lemma}[section]

\newcommand{\eqnsection}{
   \renewcommand{\theequation}{\thesection.\arabic{equation}}
   \makeatletter
   \csname @addtoreset\endcsname{equation}{section} 
   \makeatother}
    
\def \ov{\overline}

\def \be{\begin{equation}}
\def \ee{\end{equation}} 
\def \bt{\begin{theorem}} 
\def \et{\end{theorem}} 
\def \bl{\begin{lemma}}     
\def \el{\end{lemma}}
\def \bea{\begin{eqnarray}}
\def \eea{\end{eqnarray}}
\def \bas{\begin{eqnarray*}}
\def \eas{\end{eqnarray*}} 

\def \al{\alpha}  
\def \bb{\beta} 
 
\def \Ga{\Gamma}
\def \de{\delta} 
\def \De{\Delta}

\def \la{\lambda}

\def \om{\omega}
\def \Om{\Omega}

\def \si{\sigma}

\def \th{\theta}

\def \ze{\zeta}

\def \ff{\infty}
\def \wh{\widehat}
\def \wt{\widetilde}

\def \rar{\rightarrow}

\def \FF{{\cal F}}

\def \({\left(}
\def \){\right)}

\def \nn{\nonumber}

\def \bc{\begin{center} }
\def \ec{\end{center} }
\def \bs{\begin{slide} }
\def \es{\end{slide} }

\def\square{{\vcenter{\vbox{\hrule height.3pt
        \hbox{\vrule width.3pt height5pt \kern5pt
           \vrule width.3pt}
        \hrule height.3pt}}}}
\def\qed{{\hfill $\square$ \bigskip}}

 \def \Rev({\mbox{\rm Rev}(}

\eqnsection

\begin{document}
\title{Loop measures without transition probabilities}

 \author{Pat Fitzsimmons\,\, \,\, Yves Le Jan\,\, \,\, Jay Rosen \thanks{Research of   J. Rosen 
  was partially supported by  grants from the National Science Foundation and
PSC CUNY. }}
\maketitle

 \begin{abstract}   
 \end{abstract}

\section{Introduction}\label{sec-intro}

To the best of our knowledge, loop measures first appeared in the work of Symanzik on Euclidean quantum field theory, \cite{Sy}, where they are referred to as `blob measures'. They next appear in the work of Lawler and Werner, \cite{LW}. In both works, the loop measure is that associated with Brownian motion. In \cite{Le Jan1} loop measures associated with a large class of Markov processes are defined and studied. In all these cases it is assumed that the underlying Markov process has transition densities. The goal of this paper is to define and study loop measures for Markov processes without transition densities.

  Let     $X\!=\!
(\Om,  \FF_{t}, X_t,\th_{t},P^x
)$ be a transient Borel right process \cite{S}   with state space $S$, which we assume to  be locally compact   with a countable base.   We use the canonical representation of $X$ in which $\Om$ is  the set of right continuous paths paths $\om:[0,\ff)\to S_{\De}=S\cup \De$ with  $\De\notin S$, and  is such that   $\om(t)=\De$ for all $t\geq \ze=\inf \{t>0\,|\om(t)=\De\}$.   Set  $X_{t}(\om)=\om(t)$. 

Let $m$ be a Borel measure on $S$ which is finite on compact sets. We assume that with respect to   $m$, $X$ has strictly positive  potential densities $u^{\al}(x,y)$, $\al\geq 0$,  which satisfy the resolvent equations. We set $u(x,y)=u^{0}(x,y)$, and  assume that   $u(x,y)$ is excessive in $x$ for each fixed $y$.

Let $h_{z}(x)=u(x,z)$. If we assume that $u$ is finite, then the $h_{z}-$transform of $X$ is a right process on $S$, see \cite[Section 62]{S}, with probabilities $P^{x/h_{z}}$. Let $Q^{z,z}=u(z,z)P^{z/h_{z}}$. We can then define the loop measure as
\begin{equation}
\mu\(F\)=\int Q^{z,z}\({F \over \ze}\)\,dm(z),\label{0.1}
\end{equation}
for any $\FF$ measurable function $F$.
Loop measures for processes with finite potential densities but without transition densities are discussed in \cite{FR}. In the present paper we assume that the potential densities $u(x,y)$ are infinite on the diagonal, but finite off the diagonal. In this case the construction of the  measures $P^{z/h_{z}}$ given in \cite{S} breaks down. Assuming   that all points are polar, we show how to construct a family of measures $Q^{z,z}, z\in S$, which generalize  the measures $Q^{z,z}=u(z,z)P^{z/h_{z}}$ in the case of finite $u(x,y)$.

After constructing  $Q^{z,z}, z\in S$ and defining the loop measure $\mu$ using (\ref{0.1}), we show how to calculate some important moments.
We assume that
\begin{equation}
\sup_{x}\int_{K} \(u(x,y)+u(y,x)\)^{2}\,dm(y)<\ff\label{fin1}
\end{equation}
for any   compact  $K\subseteq S$. For exponentially killed Brownian motion in $R^{d}$ this means that $d\leq 3$.

\bt\label{theo-momint}
For any $k\geq 2$, and   bounded measurable functions $f_{1},\ldots, f_{k}$ with compact support
\begin{eqnarray}
&&\mu\(\prod_{j=1}^{k}\int_{0}^{\ff} f_{j}(X_{t_{j}})\,dt_{j}\)
\label{0.2}\\
&&=\sum_{\pi\in\mathcal{P}_{k}^{\odot}}\int  u(x_{1}, x_{2})u(x_{2}, x_{3})\cdots u(x_{k},x_{1}) \prod_{j=1}^{k}f_{\pi_{j}}(x_{j})\,dm(x_{j}),\nn
\end{eqnarray}
where $\mathcal{P}_{k}^{\odot}$ denotes the set of permutations of $[1,k]$ on the circle. (For example,  $(1,2,3)$,  $(3,1,2)$ and $(2,3,1)$ are considered to be one permutation   $\pi\in \mathcal{P}_{3 }^{\odot}$.)
\et

Our assumption (\ref{fin1}) will guarantee that the right hand side of (\ref{0.2}) is finite.
Note that if $k=1$ our formula would give
\begin{equation}
\mu\(\int_{0}^{\ff} f(X_{t})\,dt_{j}\)=\int u(x ,x )  f (x )\,dm(x )=\ff,\label{}
\end{equation}
for any $f\geq 0$, by our assumption that the potentials $u(x,y)$ are infinite on the diagonal.

For  $f_{1},\ldots, f_{k}$ as above consider more generally the multiple integral  
 \begin{equation}
M^{f_{1},\ldots,f_{k}}_{t}= \sum_{\pi\in \mathcal{T}_{k}^{\odot}}  \int_{0\leq r_{1}\leq\cdots\leq r_{k}\leq t}   \,f_{\pi (1)}(X_{r_{1}})\cdots  \,f_{\pi (k)}(X_{r_{k}})\,dr_1  \cdots dr_k, 
\end{equation}
where $\mathcal{T}_{k}^{\odot}$ denotes the set of   translations  $\pi$  of $[1,k]$ which are cyclic mod $k$, that is, for some $i$,
$\pi(j)=j+i,\mod k, $ for all  $j=1,\ldots,k$. In the proof of Theorem \ref{theo-momint} we first show that 
\be
\mu\(M^{f_{1},\ldots,f_{k}}_{\ff}\)
 = \int  u(x_{1}, x_{2})u(x_{2}, x_{3})\cdots u(x_{k},x_{1}) \prod_{j=1}^{k}f_{j}(x_{j})\,dm(x_{j}). \label{0.2x}
\ee
(\ref{0.2}) will then follow since 
\begin{equation}
\prod_{j=1}^{k}\int_{0}^{\ff} f_{j}(X_{t_{j}})\,dt_{j}=\sum_{\pi\in\mathcal{P}_{k}^{\odot}}M^{f_{\pi_{1}},\ldots,f_{\pi_{k}}}_{\ff}.\label{0.2xa}
\end{equation}

There is a related measure which we shall use which gives finite values even for $k=1$.  Set
\begin{equation}
\nu(F)=\int  Q^{z,z}\(F \)\,dm(z).\label{49.0p}
\end{equation}
Assume that  for any $\al>0$,  any   compact  $K\subseteq S$, and any $\wt K $ which is a compact neighborhood of $K$
\begin{equation}
\sup_{ z\in \wt K^{c}, \,\,x\in K} u^{\al}(z,x)<\ff.\label{fin2}
\end{equation}

\bt\label{theo-momnu}
For any $k\geq 1$, $\al>0$, and   bounded measurable functions $f_{1},\ldots, f_{k}$ with compact support
\begin{eqnarray}
&&\nu\(\prod_{j=1}^{k}\int_{0}^{\ff} f_{j}(X_{t_{j}})\,dt_{j}e^{-\al\ze}\)
\label{0.2a}\\
&&=\sum_{\pi\in\mathcal{P}_{k}}\int  u^{\al}(z, x_{1})u^{\al}(x_{1}, x_{2})\cdots u^{\al}(x_{k},z) \prod_{j=1}^{k}f_{\pi_{j}}(x_{j})\,dm(x_{j})\,dm(z),\nn
\end{eqnarray}
where $\mathcal{P}_{k}$ denotes the set of permutations of $[1,k]$, and both sides are finite.
\et

  We call $\mu$ the loop measure of $X$ because, when $X$ has continuous paths, $\mu$ is concentrated on the set of continuous 
 loops with a distinguished starting point (since  $Q^{x,x}$ is carried by loops starting at $x$).
 Moreover, in the next Theorem we show that  it is shift invariant. More precisely, 
  let  $\rho_{u}$ denote the loop rotation defined by
\[
\rho_{u}\om (s)
	=\left\{\begin{array}{ll}
	\om (s+u \mbox{ mod } \zeta(\om)), & \mbox{ if $0\leq s<\zeta(\om)$}\\
	\De, & \mbox{otherwise.}
	\end{array} \right.
	\]
Here, for two positive numbers $a,b$ we define  $a\mbox{ mod }  b= a-mb$ for the unique positive integer $m$ such that $\,0\leq a-mb<b $. 

For the next Theorem we need an additional assumption: for any $\de>0$
and compact $K\subseteq S$
\begin{equation}
\int_{K}  P_{\de} (z,dx ) u(x ,z)\,dm(z)<\ff.\label{0.4}
\end{equation}

\bt\label{theo-inv}
 $\mu $ is invariant under  $\rho_{u}$, for any   $u$.
\et

Note that if we have transition densities $p_{\de} (z, x ) $ then
\begin{eqnarray}
\int_{K}  P_{\de} (z,dx ) u(x ,z)\,dm(z)&=&\int \int_{K} p_{\de} (z, x )   u(x ,z)\,dm(x)\,dm(z)
\label{0.6}\\
&=&\int_{K} \(\int_{\de}^{\ff} p_{t}(x ,x)\,dt\)\,dm(x).  \nonumber
\end{eqnarray}
In our work on processes with transition densities, it was  always assumed that   $\sup_{x}\int_{\de}^{\ff}  p_{t}(x ,x)\,dt<\ff$   for any $\de>0$, which indeed gives (\ref{0.4}).

 For the next Theorem we assume that the measure $m$ is excessive. With this assumption there is always a dual process $\hat X$ (essentially uniquely determined), but in general it is a moderate Markov process. We assume that the measures $\hat U(\cdot, y)$ are absolutely continuous with respect to $m$ for each $y\in S$.

For CAF's  $L^{\nu_{1}}_{t},\ldots L^{\nu_{k}}_{t}$ with  Revuz measures $\nu_{1},\dots, \nu_{k}$, let 
\begin{equation}
A^{\nu_{1},\ldots,\nu_{k}}_{t}= \sum_{\pi\in  \mathcal{T}_{k}^{\odot}}  \int_{0\leq r_{1}\leq\cdots\leq r_{k}\leq t}   \,dL^{\nu_{\pi (1)}}_{r_{1}}\cdots  \,dL^{\nu_{\pi (k)}}_{r_{k}}. \label{2.0}
\end{equation}
 We refer to $A^{\nu_{1},\ldots,\nu_{k}}_{t}$ as a multiple CAF. 

\bt\label{theo-momCAF}
For any $k\geq 2$, and any CAF's  $L^{\nu_{1}}_{t},\ldots L^{\nu_{k}}_{t}$ with  Revuz measures $\nu_{1},\dots, \nu_{k}$,  
 \be
 \mu\(A^{\nu_{1},\ldots,\nu_{k}}_{\ff}\)
 =\int  u(x_{1}, x_{2})u(x_{2}, x_{3})\cdots u(x_{k},x_{1}) \prod_{j=1}^{k}\,d \nu_{j}(x_{j}). \label{0.3x}  
\ee
and
 \be
 \mu\(\prod_{j=1}^{k}L^{\nu_{j}}_{\ff}\)
 =\sum_{\pi\in\mathcal{P}_{k}^{\odot}}\int  u(x_{1}, x_{2})u(x_{2}, x_{3})\cdots u(x_{k},x_{1}) \prod_{j=1}^{k}\,d \nu_{\pi_{j}}(x_{j}). \label{0.3}  
\ee
\et

The finiteness of the right hand side of (\ref{0.3}) will depend on the potential densities $u(x,y)$
and the measures $\nu_{1},\dots, \nu_{k}$. For a more thorough discussion see \cite[(1.5)]{LMR}
and the  paragraph there following (1.5).

With the results of this paper, most of the results of \cite{FR, LMR, LMR2} on loop measures, loop soups, CAF's and intersection local times will carry over to processes without transition densities.

\section{Construction of $Q^{z,z}$}\label{sec-const}

Let us fix  $z\in S$ and consider the excessive function $h_{z}(x):=u(x,z)$, finite and strictly positive on the subspace $S_z:=\{x\in S: x\not=z\}$.
Doob's $h$-transform theory yields the existence of  laws  $P^{x,z}$, $x\in S_z$, on path space under which the coordinate process is Markov with transition semigroup
\be
P_t^z(x,dy) := P_t(x,dy){h_{z}(y)\over h_{z}(x)}.\label{q1}
\ee
See, for example,  \cite[pp.~298--299]{S}. Now consider the family of measures
\be
\eta^z_t(dx) := P_t(z,dx)h_{z}(x).\label{q2}
\ee
Since we assume that the singleton $\{z\}$ is polar, the transition semigroup $(P_t)$ will not charge $\{z\}$,
so these may be viewed as measures on $S_z$ or on $S$. Adopting the latter point of view, it is immediate that $(\eta^z_t)_{t>0}$ is an entrance law 
for $(P^z_t)$. There is a general theorem guaranteeing the existence of a right process with one-dimensional distributions $(\eta^z_t)$ and transition semigroup
$(P^z_t)$; see  \cite[Proposition (3.5)]{GG}.
The law of this process is the desired $Q^{z,z}$. Aside from the entrance law identity $\eta^z_tP_s^z =\eta_{t+s}^z$, their result only requires that each of the measures $\eta^z_t$ be $\sigma$-finite, which is clearly the case in the present discussion. 

  With this we immediately obtain, for $0<t_{1}<\cdots<t_{k}$, 
\begin{eqnarray}
Q^{z,z}\(\prod_{j=1}^{k}f_{j}(X_{t_{j}})\)&=&\eta^z_{t_{1}}(dx_{1})\prod_{j=2}^{k} P_{t_{j}-t_{j-1}}^z(x_{j-1},dx_{j})\prod_{j=1}^{k}f_{j}(x_{j})
\label{49.1}\\
&=& \prod_{j=1}^{k} P_{t_{j}-t_{j-1}} (x_{j-1},dx_{j})\prod_{j=1}^{k}f_{j}(x_{j}) u(x_{k},z) \nonumber
\end{eqnarray}
 with $t_{0}=0$ and $x_{0}=z$. Hence 
\begin{eqnarray}
&&
Q^{z,z}\(\int_{0<t_{1}<\cdots<t_{k}<\ff}\prod_{j=1}^{k} f_{j}(X_{t_{j}})\,dt_{j}\)\label{49.2}\\
&&=\int u(z,x_{1})u(x_{1}, x_{2})\cdots u(x_{k},z) \prod_{j=1}^{k}f_{j}(x_{j})\,dm(x_{j}),\nn
\end{eqnarray}
so that
\begin{eqnarray}
&&
Q^{z,z}\(\prod_{j=1}^{k}\int_{0}^{\ff} f_{j}(X_{t_{j}})\,dt_{j}\)\label{49.3}\\
&&=\sum_{\pi\in\mathcal{P}_{k}}\int u(z,x_{1})u(x_{1}, x_{2})\cdots u(x_{k},z) \prod_{j=1}^{k}f_{\pi_{j}}(x_{j})\,dm(x_{j}).\nn
\end{eqnarray}

Returning to (\ref{49.1}) with $0<t_{1}<\cdots<t_{k}$ and using the fact that $\ze>t_{k}$ implies that $\ze=t_{k}+\ze\circ \th_{t_{k}}$ we have 
\begin{eqnarray}
&&
Q^{z,z}\(\prod_{j=1}^{k}f_{j}(X_{t_{j}})e^{-\al \ze}\) \label{49.4}\\
&&=  Q^{z,z}\(\prod_{j=1}^{k}f_{j}(X_{t_{j}})e^{-\al t_{k}}\(e^{-\al \ze}\circ \th_{t_{k}}\)\)
\nn\\
&&= \prod_{j=1}^{k}  P^{\al}_{t_{j}-t_{j-1}} (x_{j-1},dx_{j})\prod_{j=1}^{k}f_{j}(x_{j}) h_{z}(x_{k}) P^{x_{k},z}\(e^{-\al \ze}\).\nonumber
\end{eqnarray}

Note that by (\ref{q1}) and the fact that  $X$ has $\al$-potential densities for all $\al\geq 0$
\begin{eqnarray}
P^{x_{k},z}\(\int_{0}^{\ff }e^{-\al t}1_{\{S\}}(X_{t})\,dt\)&=&\int_{0}^{\ff }e^{-\al t}P^{x_{k},z}\(1_{\{S\}}(X_{t})\)\,dt
\label{q3}\\
&=&\int_{0}^{\ff }e^{-\al t}\int_{S} P_t(x_{k},dy)    {h_{z}(y)\over h_{z}(x_{k})}\,dt\nonumber\\
&=&\int u^{\al}(x_{k},y) {h_{z}(y)\over h_{z}(x_{k})} \,dm(y)   .\nonumber
\end{eqnarray}
Combining this with  our assumption that the   $\al$-potential densities satisfy the resolvent equation we see that
\begin{eqnarray}
&& h_{z}(x_{k}) P^{x_{k},z}\(e^{-\al \ze}\)= h_{z}(x_{k}) P^{x_{k},z}\(1-\al\int_{0}^{\ze }e^{-\al t}\,dt\)
\label{49.5}\\
&& =  h_{z}(x_{k}) -\al h_{z}(x_{k}) P^{x_{k},z}\(\int_{0}^{\ff }e^{-\al t}1_{\{S\}}(X_{t})\,dt\) \nonumber\\
&& =  u(x_{k},z) -\al \int u^{\al}(x_{k},y) u(y ,z)\,dm(y)=u^{\al}(x_{k},z). \nonumber
\end{eqnarray}
Using this in (\ref{49.4}) we obtain
\begin{eqnarray}
&&
Q^{z,z}\(\prod_{j=1}^{k}f_{j}(X_{t_{j}})e^{-\al \ze}\) \label{49.4f}\\
&&= \prod_{j=1}^{k}  P^{\al}_{t_{j}-t_{j-1}} (x_{j-1},dx_{j})u^{\al}(x_{k},z)\prod_{j=1}^{k}f_{j}(x_{j}) .\nonumber
\end{eqnarray}

We then  have 
\begin{eqnarray}
&&
Q^{z,z}\(e^{-\al \ze}\int_{0<t_{1}<\cdots<t_{k}<\ff}\prod_{j=1}^{k} f_{j}(X_{t_{j}})\,dt_{j}\)\label{49.6}\\
&&=\int u^{\al}(z,x_{1})u^{\al}(x_{1}, x_{2})\cdots u^{\al}(x_{k},z) \prod_{j=1}^{k}f_{j}(x_{j})\,dm(x_{j}),\nn
\end{eqnarray} 
and consequently
\begin{eqnarray}
&&
Q^{z,z}\(e^{-\al \ze}\prod_{j=1}^{k}\int_{0}^{\ff} f_{j}(X_{t_{j}})\,dt_{j}\)\label{49.7}\\
&&=\sum_{\pi\in\mathcal{P}_{k}}\int u^{\al}(z,x_{1})u^{\al}(x_{1}, x_{2})\cdots u^{\al}(x_{k},z) \prod_{j=1}^{k}f_{\pi_{j}}(x_{j})\,dm(x_{j}).\nn
\end{eqnarray}

\section{The loop measure and its moments}\label{sec-loopmom}

Set
\begin{equation}
\mu(F)=\int  Q^{z,z}\({F \over \ze}\)\,dm(z).\label{49.0}
\end{equation}

{\bf  Proof of Theorem \ref{theo-momint}: }
We   use an argument from the proof of  \cite[Lemma 2.1]{FR}, which is due to Symanzik, \cite{Sy}.

It follows from the resolvent equation that the potential densities $u^{\bb}(x ,y)$ are continuous and monotone decreasing in $\bb$, for $x\neq y$.
  Using this together with the resolvent equation and the monotone convergence theorem we obtain that for $x_{k}\neq x_{1}$
 \begin{equation}
\int_{S}  u^{\al}(x_{k},z) u^{\al}(z,x_{1})\,dm (z)=-  {d \over d\al} u^{\al}(x_{k}, x_{1}).  \label{49.9}
\ee
Hence using  (\ref{49.6})
\begin{eqnarray}
&&\hspace{-.3 in}
\int Q^{z,z}\(e^{-\al \ze}M^{f_{1},\ldots f_{k}}_{\ff}\)\,dm(z)\label{49.10}\\
&&\hspace{-.3 in}=-\sum_{\pi\in\mathcal{T}^{\odot}_{k}}\int  u^{\al}(x_{1}, x_{2}) u^{\al}(x_{2}, x_{3})\cdots u^{\al}(x_{k-1}, x_{k}){d \over d\al}u^{\al}(x_{k}, x_{1}) \prod_{j=1}^{k}f_{\pi_{j}}(x_{j})\,dm(x_{j})\nn\\
&&\hspace{-.3 in}=-{d \over d\al} \int  u^{\al}(x_{1}, x_{2}) u^{\al}(x_{2}, x_{3})\cdots u^{\al}(x_{k-1}, x_{k})u^{\al}(x_{k}, x_{1}) \prod_{j=1}^{k}f_{j}(x_{j})\,dm(x_{j}).\nn
\end{eqnarray}
For the last step we used the product rule for differentiation and the fact that in the middle line we are summing over all translations mod k.

Since, as mentioned, $u^{\al}(x, y)$ is monotone decreasing in $\al$ for $x\neq y$,
\begin{equation}
v(x, y)= \lim_{\al\rar\ff} \,\,u^{\al}(x, y)\label{}
\end{equation}
exists and 
\begin{equation}
\int v(x, y)f(y)\,dm(y)=\lim_{\al\rar\ff}\int_{0}^{\ff }e^{-\al t}\int  P_t(x_{k},dy)  f(y)\,dt=0.\label{}
\end{equation}
Hence $v(x, y)=0$ for $m-$a.e. $y$.
Integrating (\ref{49.10}) with respect to $\al$ from $0$ to $\ff$ and using Fubini's theorem 
 we then obtain (\ref{0.2x}).
(\ref{0.2})   then follows by (\ref{0.2xa}).

To show that the right hand side of (\ref{0.2}) is finite  we repeatedly use the Cauchy-Schwarz inequality and our assumption (\ref{fin1}). See  the proof of \cite[Lemma 3.3]{LMR}.\qed

{\bf  Proof of Theorem \ref{theo-momnu}: } The formula (\ref{0.2a}) follows immediately from  (\ref{49.7}). When $k\geq 2$, the right hand side of (\ref{0.2a}) can be shown to be finite by repeatedly using the Cauchy-Schwarz inequality, our assumption (\ref{fin1}) and the fact that $u^{\al}( x, z)$ is integrable in $z$ for any $\al>0$. When $k=1$, if $K$ is a compact set containing the support of $f_{1}$ and $\wt K $ is a compact neighborhood of $K$, then
\bea
&&\int \int u^{\al}(z,x)u^{\al}(x,z)f_{1}(x)\,dm(x)\,dm(z)\nn\\
&&=\int_{\wt K} \int u^{\al}(z,x)u^{\al}(x,z)f_{1}(x)\,dm(x)\,dm(z)\nn\\
&&+\int_{\wt K^{c}}  \int u^{\al}(z,x)u^{\al}(x,z)f_{1}(x)\,dm(x)\,dm(z).\nn
\eea
Using  (\ref{fin1})
\begin{eqnarray}
&&\int_{\wt K} \(\int u^{\al}(z,x)u^{\al}(x,z)f_{1}(x)\,dm(x)\)\,dm(z)
\label{}\\
&& \leq m(\wt K) \sup_{z} \int u^{\al}(z,x)u^{\al}(x,z)f_{1}(x)\,dm(x)<\ff,\nonumber
\end{eqnarray}
  and using   (\ref{fin2})
\begin{eqnarray}
&&\int_{\wt K^{c}}  \int u^{\al}(z,x)u^{\al}(x,z)f_{1}(x)\,dm(x)\,dm(z)
\label{}\\
&&\leq C   \int  \(\int u^{\al}(x,z)\,dm(z)\) f_{1}(x)\,dm(x)<\ff.\nonumber
\end{eqnarray}
\qed

\section{Subordination}\label{sec-sub}

The basic idea in our proof that the loop measure is shift invariant is to show that the loop measure can be obtained as the `limit' of loop measures for processes with transition densities. These processes will be obtained from the original process by subordination. 

  We consider a subordinator $T_{t}$ which is a compound Poisson process with Levy measure $c\psi$ so that
\be 
E^{x}\(f\(X_{T_{t}}\)\)=\sum_{j=1}^{\ff}{(ct)^{j} \over j!}e^{-ct}\int_{0}^{\ff} E^{x}\(f\(X_{s}\)\)\psi^{*j}(ds).
\label{10.1}
\ee
If we take $\psi$ to be exponential with parameter $\th$, then $\psi^{*j}(ds)={s^{j-1}\th^{j} \over \Ga(j)}e^{-s\th}\,ds$ so that we have
\be 
E^{x}\(f\(X_{T_{t}}\)\)=\sum_{j=1}^{\ff}{(ct)^{j} \over j!}e^{-ct}\int_{0}^{\ff} E^{x}\(f\(X_{s}\)\){s^{j-1}\th^{j} \over \Ga(j)}e^{-s\th}\,ds.
\label{10.2}
\ee
Hence the subordinated transition semigroup
\be 
\wt P_{t}(x,\,dy)=\sum_{j=1}^{\ff}{(ct)^{j} \over j!}e^{-ct}\int_{0}^{\ff} P_{s}(x,\,dy){s^{j-1}\th^{j} \over \Ga(j)}e^{-s\th}\,ds.
\label{10.3}
\ee
Noting that 
\begin{equation}
\hspace{-.1 in}\int_{0}^{\ff} P_{s}(x,A) s^{j-1} e^{-s\th}\,ds={d^{j-1} \over d\th^{j-1} }\int_{0}^{\ff} P_{s}(x,A) e^{-s\th}\,ds={d^{j-1} \over d\th^{j-1} }U^{\th}(x,A),\label{10.4}
\end{equation}
we see that $\wt P_{t}(x,\,dy)$ is absolutely continuous with respect to the measure $m$ on $S$, and we can choose transition densities $\wt p_{t}(x,y)$.

From now on we take $\th=c=n$, and use $(n)$ as a superscript or subscript to denote objects with respect to the subordinated process, denoted by  $    X^{(n)}_{t}$.

\bl\label{lem-potn}
\begin{equation}
    u_{(n)}^{\al}(x,y)={1 \over (1+\al/n)^{2}} u^{ \al/ (1+\al/n)}(x,y).\label{10.6}
\end{equation}
\el
In particular, 
\begin{equation}
    u_{(n)}^{\al}(x,y)\leq     u_{(n)}(x,y)=u(x,y).\label{10.6a}
\end{equation} 

{\bf  Proof: }
\bea 
    U_{(n)}^{\al}(x,\,dy)&=&
\int_{0}^{\ff}e^{-\al t}     P^{(n)}_{t}(x,\,dy)\,dt\label{10.6r}\\
&=&{1 \over \al+n}\sum_{j=1}^{\ff}\int_{0}^{\ff} \({n \over \al+n}\)^{j}P_{s}(x,\,dy){s^{j-1}n^{j} \over \Ga(j)}e^{-sn}\,ds.
\nn\\
&=&{n^{2} \over (\al+n)^{2}}\int_{0}^{\ff}P_{s}(x,\,dy)\(\sum_{j=1}^{\ff} {s^{j-1}(n^{2}/(\al+n))^{(j-1)} \over \Ga(j)}e^{-sn}\)\,ds
\nn\\
&=&{n^{2} \over (\al+n)^{2}}\int_{0}^{\ff}\(e^{sn^{2}/(\al+n)}e^{-sn}\)P_{s}(x,\,dy)\,ds\nn\\
&=&{n^{2} \over (\al+n)^{2}}\int_{0}^{\ff}e^{-s\(n-{n^{2} \over \al+n}\)} P_{s}(x,\,dy)\,ds\nn\\
&=&{n^{2} \over (\al+n)^{2}} u^{ n\al/ ( \al+n)}(x,y)\,dm(y)\nn\\
&=&{1 \over (1+\al/n)^{2}} u^{ \al/ (1+\al/n)}(x,y)\,dm(y).\nn
\eea
\qed

\bl\label{lem-lapn}
For any  $\al, \al_{j}, j=1,\ldots, k$ and continuous compactly supported $f_{j}, j=1,\ldots, k$
\bea
&&
\lim_{n\to\ff}\,\int_{R_{+}^{k}} \prod_{j=1}^{k}e^{-\al_{j}t_{j}}    Q_{(n)}^{z,z}\(\prod_{j=1}^{k}f_{j}(    X^{(n)}_{t_{j}})e^{-\al \ze}\)\prod_{j=1}^{k}\,dt_{j} \label{10.10}\\
&& = \int_{R_{+}^{k}}\prod_{j=1}^{k}e^{-\al_{j}t_{j}}Q^{z,z}\(\prod_{j=1}^{k}f_{j}(X_{t_{j}})e^{-\al \ze}\)\prod_{j=1}^{k}\,dt_{j}.\nn
\eea
\el

{\bf  Proof: }  Recall  (\ref{49.4f}). For $0<t_{1}<\cdots<t_{k}$
\begin{eqnarray}
&&
Q^{z,z}\(\prod_{j=1}^{k}f_{j}(X_{t_{j}})e^{-\al \ze}\) \label{10.7}\\
&&= \prod_{j=1}^{k}  P^{\al}_{t_{j}-t_{j-1}} (x_{j-1},dx_{j})u^{\al}(x_{k},z)\prod_{j=1}^{k}f_{j}(x_{j})\nonumber
\end{eqnarray}
with the corresponding
\begin{eqnarray}
&&
    Q_{(n)}^{z,z}\(\prod_{j=1}^{k}f_{j}(    X^{(n)}_{t_{j}})e^{-\al \ze}\) \label{10.8}\\
&&= \prod_{j=1}^{k}      P^{(n), \al}_{t_{j}-t_{j-1}} (x_{j-1},dx_{j})    u_{(n)}^{\al}(x_{k},z)\prod_{j=1}^{k}f_{j}(x_{j}).\nonumber
\end{eqnarray}

Using (\ref{10.8}) we see that
\begin{eqnarray}
&&\int_{R_{+}^{k}} \prod_{j=1}^{k}e^{-\al_{j}t_{j}}    Q_{(n)}^{z,z}\(\prod_{j=1}^{k}f_{j}(    X^{(n)}_{t_{j}})e^{-\al \ze}\) \prod_{j=1}^{k}\,dt_{j} 
\label{10.11}\\
&&=\sum_{\pi\in \mathcal{P}_{k}}\int_{\{0\leq t_{1}\leq\cdots\leq t_{k}\leq \ff\}} \prod_{j=1}^{k}e^{-\al_{j}t_{j}}    Q_{(n)}^{z,z}\(\prod_{j=1}^{k}f_{\pi(j)}(    X^{(n)}_{t_{j}})e^{-\al \ze}\)\prod_{j=1}^{k}\,dt_{j} 
\nn\\
&& =\sum_{\pi\in \mathcal{P}_{k}}\int_{\{0\leq t_{1}\leq\cdots\leq t_{k}\leq \ff\}} \prod_{j=1}^{k}      P^{(n),\al+\sum_{l=j}^{k}\al_{l} }_{t_{j}-t_{j-1}} (x_{j-1},dx_{j})\nn\\
&&\hspace{2 in}    u_{(n)}^{\al}(x_{k},z)\prod_{j=1}^{k}f_{\pi(j)}(x_{j}) \prod_{j=1}^{k}\,dt_{j} \nonumber\\
&&=\sum_{\pi\in \mathcal{P}_{k}}  \prod_{j=1}^{k}    U_{(n)}^{\al+\sum_{l=j}^{k}\al_{l} }  (x_{j-1},dx_{j})     u_{(n)}^{\al}(x_{k},z)\prod_{j=1}^{k}f_{\pi(j)}(x_{j})
\nn \\
&&=\sum_{\pi\in \mathcal{P}_{k}} \int  \prod_{j=1}^{k}    u_{(n)}^{\al+\sum_{l=j}^{k}\al_{l} }  (x_{j-1},x_{j})     u_{(n)}^{\al}(x_{k},z)\prod_{j=1}^{k}f_{\pi(j)}(x_{j})\,dm(x_{j}).
\nn
\end{eqnarray}
  (\ref{10.10})  now follows from (\ref{10.6a}), (\ref{fin1}) and the dominated convergence theorem.
\qed

It follows from (\ref{10.10}) that for a.e. $t_{1},\cdots,t_{k}$
\begin{equation}
\lim_{n\to\ff}\,    Q_{(n)}^{z,z}\(\prod_{j=1}^{k}f_{j}(    X^{(n)}_{t_{j}})e^{-\al \ze}\) =Q^{z,z}\(\prod_{j=1}^{k}f_{j}(X_{t_{j}})e^{-\al \ze}\). \label{10.9}
\end{equation}

\bl\label{lem-lapnu} 
For any  $\al, \al_{j}, j=1,\ldots, k$ and continuous compactly supported $f_{j}, j=1,\ldots, k$
\bea
&&
\lim_{n\to\ff}\,\int_{R_{+}^{k}} \prod_{j=1}^{k}e^{-\al_{j}t_{j}}   \nu_{(n)} \(\prod_{j=1}^{k}f_{j}(    X^{(n)}_{t_{j}})e^{-\al \ze}\)\prod_{j=1}^{k}\,dt_{j} \label{10.10a}\\
&& = \int_{R_{+}^{k}}\prod_{j=1}^{k}e^{-\al_{j}t_{j}}\nu \(\prod_{j=1}^{k}f_{j}(X_{t_{j}})e^{-\al \ze}\)\prod_{j=1}^{k}\,dt_{j}.\nn
\eea
\el

{\bf  Proof: }By (\ref{10.11})
\begin{eqnarray}
&&\int_{R_{+}^{k}} \prod_{j=1}^{k}e^{-\al_{j}t_{j}}   \nu_{(n)} \(\prod_{j=1}^{k}f_{j}(    X^{(n)}_{t_{j}})e^{-\al \ze}\)\prod_{j=1}^{k}\,dt_{j}
\label{10.11c}\\
&&=\sum_{\pi\in \mathcal{P}_{k}} \int  \prod_{j=1}^{k}    u_{(n)}^{\al+\sum_{l=j}^{k}\al_{l} }  (x_{j-1},x_{j})     u_{(n)}^{\al}(x_{k},z)\prod_{j=1}^{k}f_{\pi(j)}(x_{j})\,dm(x_{j})\,dm(z).   \nonumber
\end{eqnarray}
If $k\geq 2$, then using the resolvent equation we see that
\begin{eqnarray}
&&\(\sum_{l=1}^{k}\al_{l}\)\int_{R_{+}^{k}} \prod_{j=1}^{k}e^{-\al_{j}t_{j}}   \nu_{(n)} \(\prod_{j=1}^{k}f_{j}(    X^{(n)}_{t_{j}})e^{-\al \ze}\)\prod_{j=1}^{k}\,dt_{j}
\label{10.11d}\\
&&=\sum_{\pi\in \mathcal{P}_{k}} \int  \prod_{j=2}^{k}    u_{(n)}^{\al+\sum_{l=j}^{k}\al_{l} }  (x_{j-1},x_{j})     u_{(n)}^{\al}(x_{k},x_{1})\prod_{j=1}^{k}f_{\pi(j)}(x_{j})\,dm(x_{j})    \nonumber\\
&&-\sum_{\pi\in \mathcal{P}_{k}} \int  \prod_{j=2}^{k}    u_{(n)}^{\al+\sum_{l=j}^{k}\al_{l} }  (x_{j-1},x_{j})     u_{(n)}^{\al+\sum_{l=1}^{k}\al_{l}}(x_{k},x_{1})\prod_{j=1}^{k}f_{\pi(j)}(x_{j})\,dm(x_{j}),    \nonumber
\end{eqnarray}
and (\ref{10.10a}) for $k\geq 2$ then follows from (\ref{10.6a}), (\ref{fin1})  and the dominated convergence theorem. Here we repeatedly use the Cauchy-Schwarz inequality. See  the proof of \cite[Lemma 3.3]{LMR}.

When $k=1$, 
\bea
&&
\int_{R_{+}}e^{-\al_{1}t_{1}}   \nu_{(n)} \(f_{1}(    X^{(n)}_{t_{1}})e^{-\al \ze}\) \,dt_{1}\label{prob1}\\
&&=\int \int    u_{(n)}^{\al+\al_{1}}(z,x)    u_{(n)}^{\al }(x,z)f_{1}(x)\,dm(x)\,dm(z).\nn
\eea
Note that by (\ref{10.6}) we have $ u_{(n)}^{\al }(x,z)\leq  u^{\al/2 }(x,z)$ for $n$ sufficiently large.
We can then use  the argument from the  proof of Theorem \ref{theo-momnu}   and the dominated convergence theorem to get (\ref{10.10a}) for $k=1$.
\qed

It follows from (\ref{10.10a}) that
\bea
&&
\lim_{n\to\ff}\,\int_{R_{+}^{k}} \prod_{j=1}^{k}e^{-\al_{j}t_{j}}   \nu_{(n)} \(\prod_{j=1}^{k}f_{j}(    X^{(n)}_{t_{j}})g(\ze)\)\prod_{j=1}^{k}\,dt_{j} \label{10.10b}\\
&& = \int_{R_{+}^{k}}\prod_{j=1}^{k}e^{-\al_{j}t_{j}}\nu \(\prod_{j=1}^{k}f_{j}(X_{t_{j}})g(\ze)\)\prod_{j=1}^{k}\,dt_{j}\nn
\eea
for all continuous exponentially bounded functions $g$. We will be particularly interested in $g$ of the form
\begin{equation}
g(\ze)={\prod_{j=1}^{k} (1-e^{-\bb_{j} \zeta})\over \prod_{j=1}^{k}(1-e^{-\alpha_{j} \zeta})}e^{-\al \ze}h_{s}(\ze)\label{gform}
\end{equation}
where $0\leq h_{s}(\ze)\leq 1$ is a continuous function with $h_{s}(\ze)=0$ for $\ze\leq s$.

\section{Invariance under loop rotation}\label{sec-inv}

{\bf  Proof of Theorem \ref{theo-inv}: }
Because the lifetime $\zeta$ is rotation invariant ($\zeta(\rho_v\omega) =\zeta(\omega)$ so long as $\zeta(\omega)<\infty$), the rotation invariance of the loop measure 
\begin{equation}
\mu(F)=\nu\({F \over \ze}\)\label{49.0pq}
\end{equation}
 is equivalent to that of the measure $\nu$.

We note that by (\ref{0.4})
\be 
\nu\( f (X_{\de})\)=\int  P_{\de} (z,dx ) f (x) u(x ,z)\,dm(z)<\ff
\label{conj2}
\ee
 for any $\de>0$ and bounded measurable $f$ with compact support.

We next recall some ideas from \cite{FR}. Let us define the process $\overline X$ to be the periodic extension of $X$; that is,
\begin{equation}
\ov X_t
	=\left\{\begin{array}{ll}
X_{t-q\zeta}, & \mbox{ if $q\zeta\le t<(q+1)\zeta$, $q=0,1,2,\ldots$}\\
	 X_t, & \mbox{if $\zeta=\infty$}
	\end{array} \right.
	\label{p2}
\end{equation} 

It will be convenient to write
\be
\ov I_\alpha(f):=\int_0^\infty e^{-\alpha t} f(\ov X_t)\,dt,\qquad   I_\alpha(f):=\int_0^\infty e^{-\alpha t} f(X_t)\,dt.
\label{p3}
\ee
The key observation is that
\be
\ov I_\alpha(f) ={I_\alpha(f)\over 1-e^{-\alpha \zeta}},
\label{p4}
\ee
for all $\alpha>0$. This follows from
\begin{eqnarray}
\ov I_\alpha(f):&=&\int_0^\infty e^{-\alpha t} f(\ov X_t)\,dt
\nn\\
&=& \sum_{q=0}^{\ff} \int_{q\ze}^{(q+1)\ze} e^{-\alpha t} f(\ov X_t)\,dt \nonumber\\
&=& \sum_{q=0}^{\ff}e^{-\alpha q\ze}  \int_{0}^{\ze} e^{-\alpha t} f( X_t)\,dt ={I_\alpha(f)\over 1-e^{-\alpha \zeta}}\,.\nonumber
\end{eqnarray}

The rotation invariance of $\mu$ or $\nu$ is equivalent to the following Lemma.
\bl\label{lem-rotinv}
\be
\nu\left(\prod_{j=1}^k f_j(\ov X_{t_j+r})1_{\{t_{k}<\ze\}}\right) =
\nu\left(\prod_{j=1}^k f_j(\ov X_{t_j})1_{\{t_{k}<\ze\}}\right)
 \label{p18glm}
\ee 
for all  $0<t_1<\cdots< t_k$ and $r>0$ and all $f_{j}\geq 0$ continuous with compact support.
\el

Let $0\leq h_{s}(\ze)\leq 1$ be a continuous function with $h_{s}(\ze)=0$ for $\ze\leq s$. To prove Lemma \ref{lem-rotinv} we first prove the following.

\bl\label{lem-rotinvs}  For  all $k\geq 1$,
 and  $0\leq t_1,\cdots,t_k<s$   and all $f_{j}\geq 0$ continuous with compact support
\be
\nu\left(\prod_{j=1}^k f_j(\ov X_{t_j+r})h_{s}(\ze)\right) =
\nu\left(\prod_{j=1}^k f_j(\ov X_{t_j})h_{s}(\ze)\right).
 \label{corer}
\ee 
\el

{\bf  Proof of Lemma \ref{lem-rotinvs}: } Using first (\ref{p4}) and then (\ref{10.10b}) we have that
\begin{eqnarray}
&&\lim_{n\to\ff}\,\int_{R_{+}^{k}} \prod_{j=1}^{k}e^{-\al_{j}t_{j}}   \nu^{(n)}\(\prod_{j=1}^{k}f_{j}\(\,\ov {    X^{(n)}}_{t_{j}}\)(1-e^{-\bb_{j}  \zeta})e^{-\al \ze}h_{s}(\ze)\) \prod_{j=1}^{k}\,dt_{j}
\nn\\
&&=  \lim_{n\to\ff}\,\int_{R_{+}^{k}} \prod_{j=1}^{k}e^{-\al_{j}t_{j}}   \nu^{(n)}\(\prod_{j=1}^{k}f_{j}\(\, {    X^{(n)}}_{t_{j}}\){(1-e^{-\bb_{j}  \zeta}) \over (1-e^{-\al_{j} \ze})}e^{-\al \ze}h_{s}(\ze)\) \prod_{j=1}^{k}\,dt_{j} \nn\\
&& = \int_{R_{+}^{k}}\prod_{j=1}^{k}e^{-\al_{j}t_{j}}\nu\(\prod_{j=1}^{k}f_{j}(  X_{t_{j}}){(1-e^{-\bb_{j}  \zeta}) \over(1-e^{-\al_{j} \ze})}e^{-\al \ze}h_{s}(\ze)\)\prod_{j=1}^{k}\,dt_{j}\nn\\
&& = \int_{R_{+}^{k}}\prod_{j=1}^{k}e^{-\al_{j}t_{j}}\nu\(\prod_{j=1}^{k}f_{j}( \ov X_{t_{j}})(1-e^{-\bb_{j}  \zeta}) e^{-\al \ze}h_{s}(\ze)\)\prod_{j=1}^{k}\,dt_{j}.\label{10.21h}
\eea
It follows from this that for a.e. $0\leq t_1,\cdots,t_k$
\bea
&& 
\lim_{n\to\ff}\,   \nu^{(n)}\(\prod_{j=1}^{k}f_{j}\(\,\ov {    X^{(n)}}_{t_{j}}\)(1-e^{-\bb_{j}  \zeta})e^{-\al \ze}h_{s}(\ze)\)\label{10.21}\\
&&\hspace{1 in} =\nu\left(\prod_{j=1}^k f_j(\ov X_{t_j})(1-e^{-\bb_{j}  \zeta})e^{-\al \ze}h_{s}(\ze)\right). \nn
\eea
The same calculations show that for any $r>0$, 
\bea
&& 
\lim_{n\to\ff}\,   \nu^{(n)}\(\prod_{j=1}^{k}f_{j}\(\,\ov {    X^{(n)}}_{t_{j}+r}\)(1-e^{-\bb_{j}  \zeta})e^{-\al \ze}h_{s}(\ze)\) \label{10.21j}\\
&&\hspace{1 in} =\nu\left(\prod_{j=1}^k f_j(\ov X_{t_j+r})(1-e^{-\bb_{j}  \zeta})e^{-\al \ze}h_{s}(\ze)\right) \nn
\eea
for a.e. $0\leq t_1,\cdots,t_k$.
Since $   \nu^{(n)}$ is invariant under loop rotation, see \cite[Lemma 2.4]{FR} for the simple proof, it follows from our last two displays that for a.e. $0\leq t_1,\cdots,t_k$
\bea
&&
\nu\left(\prod_{j=1}^k f_j(\ov X_{t_j+r})(1-e^{-\bb_{j}  \zeta})e^{-\al \ze}h_{s}(\ze)\right) \label{corer1}\\
&&  \hspace{1 in} =
\nu\left(\prod_{j=1}^k f_j(\ov X_{t_j})(1-e^{-\bb_{j}  \zeta})e^{-\al \ze}h_{s}(\ze)\right).
\nn
\eea 

We now  use an argument from \cite{FR} (see from (5.31) there until the end of the paragraph). 
By Fubini we can find a set $T\subseteq (0,s)$ with full  measure   such that for all $t_{1}\in T$ we have that (\ref{corer1}) holds for a.e.  $t_2,\ldots,t_k\in (0,s)$. Using (\ref{conj2}) with $\de=t_{1}$, the boundedness and  continuity of the $f_{j}$ and the right continuity of $\bar X_{t}$ it follows from the Dominated Convergence Theorem that  (\ref{corer1}) holds for all  $(t_{1},t_2,\ldots,t_k)\in T\times [0,s)^{k-1}$. Let now $f_{1,n}$ be a sequence of continuous functions with compact support with the property that $f_{1,n}\uparrow 1$. By the above, (\ref{corer1}) with $f_{1}$ replaced by $f_{1,n}$ holds for all  $(t_{1},t_2,\ldots,t_k)\in T_{n}\times [0,s)^{k-1}$ for an appropriate $T_{n}\subseteq (0,s)$ with full  measure . In particular $T_{\ast}=\cap_{n}T_{n}\not=\emptyset$, and we can apply the Monotone Convergence Theorem with $t_{1}\in T_{\ast}$ to conclude  that 
\bea
&&
\nu\left((1-e^{-\bb_{1}  \zeta})\prod_{j=2}^k f_j(\ov X_{t_j+r})(1-e^{-\bb_{j}  \zeta})e^{-\al \ze}h_{s}(\ze)\right)\label{res.17}\\&&=\nu\left((1-e^{-\bb_{1}  \zeta})\prod_{j=2}^k f_j(\ov X_{t_j})(1-e^{-\bb_{j}  \zeta})e^{-\al \ze}h_{s}(\ze)\right)\nn
\eea
for all $t_2,\ldots,t_k<s$. Applying once again the Monotone Convergence Theorem for $\bb_{j}\to\ff, \al\to 0$ we obtain
\begin{equation}
\nu\left(\prod_{j=2}^k f_j(\ov X_{t_j+r})h_{s}(\ze)\right)=\nu\left(\prod_{j=2}^k f_j(\ov X_{t_j})h_{s}(\ze)\right)\label{res.17d}
\end{equation} 
for all $t_2,\ldots,t_k<s$. Since $k$ is arbitrary, we obtain our Lemma.\qed

{\bf  Proof of Lemma \ref{lem-rotinv}: } Fix $0<t_1<\cdots< t_k$. Choose a sequence $s_{n}\downarrow t_{k}$. It is clear that we can choose $h_{s_{n}}$ so that $h_{s_{n}}(\ze)\uparrow 1_{\{t_{k}<\ze\}}$. Lemma \ref{lem-rotinv} then follows from Lemma \ref{lem-rotinvs} by the 
Monotone Convergence Theorem.\qed

\section{The loop measure and continuous additive functionals}\label{sec-CAF}

Before proving Theorem \ref{theo-momCAF} we will need two facts about continuous additive functionals (CAFs). The first says that to each CAF  $A$ of $X$ is associated a measure $\nu_{A}$
on $S$ such that for any measurable function $f$
\begin{equation}
U_{A}f(x):=E^x\int_0^\infty f(X_t)\,dA_t =\int_S u(x,y)f(y)\,\nu_{A}(dy),\qquad\forall x\in S.
\label{49.22}
\end{equation}
$\nu_{A}$ is referred to as the Revuz measure of $A$. The second fact we need is that if   a CAF has Revuz measure $\nu$ with respect to $X$, it has Revuz measure $h \cdot\nu$ with respect to the $h$-transform of  $X$. Following the proof of Theorem \ref{theo-momCAF} we will discuss these facts and provide references.


{\bf  Proof of Theorem \ref{theo-momCAF}: } To prove (\ref{0.3}) it is enough to prove the additive functional version of (\ref{49.2}). 
We consider first our Borel right process $X$. These considerations will then be applied to the $h$-transform of $X$ using $h_{z}=u(\cdot,z)$ for fixed $z\in S$.

Let $A^j$ ($j=1,2,\ldots$) be CAFs of $X$ with Revuz measures $\nu_j$.
Using the Markov property, see for example Theorems 28.7 and 22.8 of \cite{S}, and (\ref{49.22}) at the last step
\bea
&&
E^x\int_{\{0<t_1<t_2<\cdots<t_n<\infty\}}\prod_{j=1}^n dA^j_{t_j}\label{49.23}\\
&&
=E^x\(\int_{0}^{\ff }\(\int_{\{0<t_2<\cdots<t_{n}<\infty\}}\prod_{j=2}^{n} dA^j_{t_j}\)\circ\th_{t_{1}}\,dA^1_{t_1}\)
\nn\\
&&
=E^x\(\int_{0}^{\ff }E^{X_{t_{1}}}\(\int_{\{0<t_2<\cdots<t_{n}<\infty\}}\prod_{j=2}^{n} dA^j_{t_j}\)\,dA^1_{t_1}\)
\nn\\
&&
=\int_{S} 
u(x,x_1)E^{x_1}\(\int_{\{0<t_2<\cdots<t_{n}<\infty\}}\prod_{j=2}^{n} dA^j_{t_j}\) \nu_{1}(dx_1),
\nn
\eea
and then by induction
\bea
&&
E^x\int_{\{0<t_1<t_2<\cdots<t_n<\infty\}}\prod_{j=1}^n dA^j_{t_j}\label{49.24}\\
&&
=\int_{S^n} u(x,x_1)u(x_1,x_2)\cdots u(x_{n-1},x_n)\prod_{j=1}^n\nu_j(dx_j).
\nn
\eea
Notice that  by our assumption that $u(x,x_1)$ is excessive in $x$ for each $x_{1}$, the expressions in (\ref{49.24}) are excessive functions of $x$. Thus if $\eta = (\eta_t)$ is an entrance law, then
writing $E^\eta$ for the measure under which the one-dimensional distributions are given by the entrance law we have
 \begin{eqnarray}
&&
E^\eta \int_{\{0<t_1<t_2<\cdots<t_n<\infty\}}\prod_{j=1}^n dA^j_{t_j}\label{49.24a}\\
&&=\uparrow\lim_{t\downarrow 0}\int_S\eta_t(dx) E^x\int_{\{0<t_1<t_2<\cdots<t_n<\infty\}}\prod_{j=1}^n dA^j_{t_j} \nn\\
&&=\uparrow\lim_{t\downarrow 0}\int_S\eta_t(dx)\int_{S^n} u(x,x_1)u(x_1,x_2)\cdots u(x_{n-1},x_n)\prod_{j=1}^n\nu_j(dx_j)  \nn\\
&&=\int_{S^n} g(x_1)u(x_1,x_2)\cdots u(x_{n-1},x_n)\prod_{j=1}^n\nu_j(dx_j), 
\nn
\end{eqnarray}
where $g(x_1):=\uparrow\lim_{t\downarrow 0}\int\eta_t(dx)u(x,x_1)$.

Now apply the above to the $h$-transform of the original process $X$, with $h_{z}=u(\cdot,z)$ for a fixed $z\in S$, as described in Section \ref{sec-const}.
This process  has potential density $u^{h_{z}}(x,y) = u(x,y)/h_{z}(x)$ with respect to the 
 measure $h_{z}(y)\,m(dy)$. Also, if a CAF has Revuz measure $\nu$ with respect to $X$, it has Revuz measure $h_{z}\cdot\nu$ with respect to the $h$-transform process. Thus by (\ref{49.24})
 \bea
&&
E^{x,z}\int_{\{0<t_1<t_2<\cdots<t_n<\infty\}}\prod_{j=1}^n dA^j_{t_j}\label{49.24aa}\\
&&
=\int_{S^n} {u(x,x_1) \over h_{z}(x)}{u(x_1,x_2) \over h_{z}(x_{1})}\cdots {u(x_{n-1},x_n) \over h_{z}(x_{n-1})}\prod_{j=1}^n h_{z}(x_{j})\,\nu_j(dx_j)\nn\\
&&
={1 \over h_{z}(x)}\int_{S^n} u(x,x_1)u(x_1,x_2)\cdots u(x_{n-1},x_n)\,h_{z}(x_{n})\prod_{j=1}^n \,\nu_j(dx_j)
\nn\\
&&
={1 \over h_{z}(x)}\int_{S^n} u(x,x_1)u(x_1,x_2)\cdots u(x_{n-1},x_n)\,u(x_{n},z)\prod_{j=1}^n \,\nu_j(dx_j).
\nn
\eea
 When we use the entrance law $\eta^z_t(dx) = P_t(z,dx)h_{z}(x)$, the function $g$ of the preceding paragraph is 
\begin{equation}
\uparrow\lim_{t\downarrow 0}\eta^z_t(dx)u^{h_{z}}(x,x_1)
=\uparrow\lim_{t\downarrow 0}P_t(z,dx)u(x,x_1)= u(z,x_1).
\label{49.26}
\end{equation}
Thus, using the definition of $Q^{z,z}$ from Section \ref{sec-const}, 
\begin{equation}
Q^{z,z}\int_{0<t_1<\cdots<t_n<\infty}\prod_{j=1}^n dA^t_{t_j}
=\int_{S^n} u(z,x_1)u(x_1,x_2)\cdots u(x_n,z)\prod_{j=1}^n \nu_j(dx_j).
\label{49.27}
\end{equation}
Similar considerations work for the $\alpha$-potentials, and the argument given in the proof of Theorem \ref{theo-momint} proves (\ref{0.3}).\qed

We now discuss the facts mentioned at the beginning of this section.

Given a right-continuous strong Markov process X (more precisely, a Borel right Markov process) and an excessive measure m, there is always a dual process $\hat X$ (essentially uniquely determined), but in general it is a moderate Markov process: the Markov property holds only at predictable times.

 In what follows $f$ and $g$ are non-negative Borel functions on $S$. By duality 
\be
\int_S f(x) Ug(x) m(dx) = \int_S \hat Uf(y) g(y) m(dy),
\label{pf.2}
\ee
 where  the kernel $\hat U$ is the potential kernel of the moderate Markov dual of $X$. 
 Under our assumptions it follows from \cite[VI, Theorem 1.4]{BG} that the potential density $u$ can be chosen so that $x\mapsto u(x,y)$ is excessive for 
each $y$, and $y\mapsto u(x,y)$ is co-excessive (that is, excessive with respect to the moderate Markov dual process $\hat X$) for each $x$.  
 (\ref{pf.2}) implies that
\be
\hat U f(y) =\int_S u(x,y)f(x)\,m(dx),
\label{pf.3}
\ee
for $m$-a.e. $y$. 
Since both sides of (\ref{pf.3}) are co-excessive,   they agree for all $y$.

 By   \cite[(5.13)]{FG} we have the Revuz formula
\be
\int_S f(x) U_Ag(x)\, m(dx) = \int_S \hat U f(y) g(y)\nu_{A}(dy),
\label{pf.1}
\ee
where $\nu_{A}$ is the Revuz measure of the CAF $A$ with respect to $m$. Feeding (\ref{pf.3}) into (\ref{pf.1}) and varying $f$ we find that
\be
U_Ag(x)=\int_S u(x,y) g(y) \,\nu_{A}(dy),
\label{pf.4}
\ee
first for $m$-a.e. $x$, then for all $x$ because both sides of (\ref{pf.4}) are excessive. 
This proves (\ref{49.22}). 

One subtlety: the laws $\hat P^x$ of $\hat X$ are only determined modulo a class of sets (``$m$-exceptional'') defined in \cite{FG}, see (3.4) for the definition of the term, and then Remark (5.14); but that class is not charged by $\nu$, so the exception causes no problem.  

To establish the second fact that we needed, 
let $\tau_{t}$ be the right continuous inverse of $A_{t}$, and let $f$ be a positive measurable function. 
 Using the change of variables formula, \cite[Chapter 6, (55.1)]{DM2} and then Fubini
\bea
E^{x/h}\(\int_{0}^\infty f(X_t)\,dA_t\) &=&E^{x/h}\(\int_0^\infty f(X_{\tau(u)})\,du\)\label{}\\
&=&\int_0^\infty E^{x/h}\(f(X_{\tau(u)})\)\,du.\nn
\eea
Using  \cite[(62.20)]{S} and then  Fubini we have
\bea\hspace{-.4 in}
\int_0^\infty E^{x/h}\(f(X_{\tau(u)})\)\,du &=&{1 \over h(x)}\int_0^\infty E^{x}\(f(X_{\tau(u)})h(X_{\tau(u)})\)\,du\label{}\\
&=&{1 \over h(x)}E^{x}\(\int_0^\infty f(X_{\tau(u)})h(X_{\tau(u)})\,du\).\nn
\eea
Using the change of variables  formula once again, the last two formulas show that
\begin{equation}
E^{x/h}\(\int_{0}^\infty f(X_t)\,dA_t\)={1\over h(x)}E^x\(\int_0^\infty f(X_t)h(X_t)\,dA_t\).\label{}
\end{equation}
Using (\ref{49.22})  we see that
\begin{equation}
E^{x/h}\(\int_{0}^\infty f(X_t)\,dA_t\)={1\over h(x)} \int_S u(x,y)f(y)h(y)\,\nu_{A}(dy).\label{}
\end{equation}
This shows that   if   a CAF has Revuz measure $\nu$ with respect to $X$, then it has Revuz measure $h \cdot\nu$ with respect to the $h$-transform of  $X$.  
 (Recall that we use ${ u(x,y)\over h(x)}$ for the potential densities of the $h$-transform process).

\def\noopsort#1{} \def\printfirst#1#2{#1}
\def\singleletter#1{#1}
            \def\switchargs#1#2{#2#1}
\def\bibsameauth{\leavevmode\vrule height .1ex
            depth 0pt width 2.3em\relax\,}
\makeatletter
\renewcommand{\@biblabel}[1]{\hfill#1.}\makeatother
\newcommand{\bysame}{\leavevmode\hbox to3em{\hrulefill}\,}

 \def\wh{\widehat}
\def\ol{\overline}

{\footnotesize

\noindent
\begin{tabular}{lll} &  Pat Fitzsimmons
     & \hskip20pt   Yves Le Jan\\  &  Department of Mathematics
     & \hskip20pt Equipe Probabilit\'es et Statistiques \\    &
University of California, San Diego
     & \hskip20pt     Universit\'e Paris-Sud, B\^atiment 425\\    &  
    La Jolla, CA 92093-0112 & \hskip20pt   91405 Orsay Cedex
 \\    &    U.S.A.
     & \hskip20pt France \\   &   
pfitzsim@ucsd.edu
     & \hskip20pt   yves.lejan@ math.u-psud.fr
\end{tabular}
\bigskip

 \begin{center}
\begin{tabular}{lll} & \hskip20pt   Jay Rosen\\  & \hskip20pt Department of Mathematics 
     \\    
     & \hskip20pt College of Staten Island, CUNY \\    
     & \hskip20pt  Staten Island, NY
10314 \\    
     & \hskip20pt U.S.A. \\   
     & \hskip20pt  jrosen30@optimum.net
\end{tabular}
\end{center}

}


\begin{thebibliography}{10}

\bibitem {BG}
R. Blumenthal and R. Getoor,
 {\em Markov {P}rocesses and {P}otential {T}heory},
 Academic Press,
 New York, 1968.
 

\bibitem {DM2}
C. Dellacherie,  and P.-A. Meyer,  (1982).
\newblock {\em Probabilities and Potential B}.
\newblock North Holland Publishing Company, Amsterdam.




\bibitem {FG}
P. Fitzsimmons, and R. K. Getoor,  
\newblock Homogeneous Random Measures and Strongly Supermedian Kernels of a Markov Process. 
{\em Electronic Journal of Probability}, {\bf  8}, 1-54, 2003.


\bibitem {FR} P. J. Fitzsimmons and J. Rosen,
Markovian loop soups: permanental processes and isomorphism theorems, \, {\it Electron. J. Probab.},\, Volume 19 (2014), no. 60, 1-30.
 


\bibitem {GG}
R. K. Getoor,  and J. Glover,  
\newblock Constructing Markov Processes with Random Times of Birth and Death.  
{\em Seminar on Stochastic Processes}, 1986, E. Cinlar, K.L. Chung, R. K. Getoor,  and J. Glover editors, Progress in Probability and Statistics, Volume 13, pp. 35-69, BirkhŠuser, Boston (1987). 


\bibitem{LW}
G. Lawler and W. Werner,  \newblock The Brownian loop soup, \newblock {\em PTRF} 44 (2004), 197--217.
 
\bibitem{Le Jan1} Y. Le Jan, \newblock
{\em Markov paths, loops and fields. }   \'{E}cole d'\'{E}t\'{e} de Probabilit\'{e}s de Saint-Flour XXXVIII - 2008. Lecture Notes in Mathematics 2026. (2011)
Springer-Verlag, Berlin-Heidelberg. 
 
 \bibitem{LMR} Y. Le Jan, M. B.  Marcus and J.~Rosen,\newblock
Permanental fields, loop soups and continuous additive functionals, http://arxiv.org/pdf/1209.1804.pdf, Version 1.


 \bibitem{LMR2} Y. Le Jan, M. B.  Marcus and J.~Rosen,\newblock 
  Intersection local times, loop soups and permanental Wick powers, http://arxiv.org/pdf/1308.2701.pdf 
  


   \bibitem{S} M. Sharpe, {\em General theory of Markov
processes},  Acad. Press, New York, (1988).
 

\bibitem{Sy} {K}. {S}y\-man\-zyk.\newblock
     {E}uclidean quantum field theory.
         In R.~Jost, editor,  {\sl Local Quantum Theory}. Academic  Press,
Reading, MA,
       1969.
 



\end{thebibliography}
  \end{document}